\documentstyle[11pt]{article}

\newcommand{\cD}{{\cal D}}

\newcommand{\cF}{{\cal F}}

\newcommand{\cH}{{\cal H}}

\newtheorem{lemma}{Lemma}

\newenvironment{proof}{{\bf
Proof}\quad}{{\hfill$\bullet$}\ \\ \ \par}

\newcommand{\B}[1]{\begin{#1}}
\newcommand{\E}[1]{\end{#1}}

\newcommand{\C}{\mbox{$\bf C$}}     
\newcommand{\Z}{\mbox{$\bf Z$}}     

\newcommand{\im}{{\rm im}\,}

\newcommand{\df}{\mbox{\,$\stackrel{\pp{\rm def}}{=}$}\,}

\newcommand{\by}[1]{\stackrel{#1}{\rightarrow}}
\newcommand{\longby}[1]{\stackrel{#1}{\longrightarrow}}

\newcommand{\tensor}{\otimes}

\newcommand{\ie}{{\it i.e.\/}\ }

\newcommand{\cf}{{\it cf.\/}\ }

\newcommand{\sZ}{\mbox{\scriptsize{$\Z$}}}

\newcommand{\pp}[1]{\mbox{$\scriptscriptstyle {#1}$}}

\newcommand{\limdir}[1]{{\displaystyle{\mathop{\rm
lim}_{\buildrel\longrightarrow\over{#1}}}}\,}

\title{Bloch's conjecture revisited}

\author{by {\sc L.Barbieri-Viale} and {\sc V.Srinivas}}

\date{March `95}

\begin{document}

\maketitle

\B{abstract}
Let $X$ be a non-singular projective complex surface.
We can show that Bloch's conjecture (\ie, that if $p_g=0$
then the Albanese kernel vanishes) is
equivalent to the following statement:\\[6pt]
 {\em If $p_g(X)=0$ then for any given Zariski open
$U\subset X$ and $\omega\in H^2(U,{\bf C})$ there is a
smaller Zariski open $V\subset U$ such that $$\omega\mid_V
= \omega'+\omega_{\sZ}$$ where $\omega'\in F^2H^2(V,{\bf
C})$ and $\omega_{\sZ}$ is integral.}
\E{abstract}

\vspace{1cm}

Let $X$ denote a complex algebraic
surface which is smooth and complete. Let $A_0(X)$ denote
the subgroup of the Chow group of zero cycles given
by cycles of degree zero, and let $J^2(X)$ be the Albanese variety of $X$,
which we regard as the intermediate Jacobian canonically associated with
the Hodge structure on $H^3(X,\Z (2))$. Let $\phi: A_0(X)\to
J^2(X)$ be induced by the canonical (surjective) mapping
to the Albanese variety. The kernel of $\phi$ is usually
called the ``Albanese kernel''.\\[5pt] {\bf Bloch's
conjecture\ \ } {\it If $p_g(X)=0$ then the
``Albanese kernel'' vanishes.}\\[4pt] The
goal of this note is to show the following: \B{thm} Let
$X$ be as above, and assume that $p_g(X)=0$. Then
$$\mbox{\em ``Albanese kernel''}\cong
\limdir{\mbox{$U\subset X$}} \frac{H^2(U,{\bf
C})}{F^2H^2(U,{\bf C})+H^2(U,{\bf Z}(2))}$$
where the limit is taken over all non-empty Zariski open
subsets of $X$.
\E{thm}
We then clearly have the following:
\B{cor} Let $X$ be as above, and assume that $p_g(X)=0$.
The ``Albanese kernel'' vanishes if and only if for any
given Zariski open $U\subset X$ and $\omega\in H^2(U,{\bf
C})$, there is a smaller Zariski open $V\subset U$ such that
$$\omega\mid_V = \omega'+\omega_{\sZ},$$ where $\omega'\in
F^2H^2(V,{\bf C})$ and $\omega_{\sZ}$ is integral (\ie,
$\omega_{\sZ}\in\im H^2(V,{\bf Z}(2))$).
\E{cor}
We will make use of the following result only in the
particular case when $Y$ is a curve.
\begin{lemma} Let $Y$ be any reduced variety over the
complex numbers. Then  \[\limdir{\mbox{$U\subset Y$}}
\frac{H^1(U,{\bf C})}{F^1H^1(U,{\bf C})+H^1(U,{\bf
Z}(1))}\;=0\] where the limit is taken over all non-empty
Zariski open subsets of $Y$.\end{lemma}
\begin{proof}
Since the statement concerns the generic points of a
variety, we may assume that $Y$ is irreducible, non-singular and
projective. For any open $U\subset Y$, let $S=Y-U$, and
let ${\rm Div}^0_S(Y)$ be the group of divisors on $Y$
which are supported on $S$, and are homologous to 0 on
$Y$. There is an exact sequence
$$0\to H^1(Y,{\bf Z}(1))\to H^1(U,{\bf Z}(1))\longby{
res} {\rm Div}^0_S(Y)\to 0$$
where $res$ denotes the sum of the
residue maps associated to components of $S$ which are
divisors on $Y$. This underlies an exact sequence of mixed
Hodge structures, where $F^1({\rm Div}^0_S(Y)\tensor\C) = {\rm
Div}^0_S(Y)\tensor\C$ (indeed, for any such  divisor
$D$, we have that $$F^1H^2_D(Y,\C (1)) \cong
F^0H^0(\tilde D,\C (0))$$ for a desingularization $\tilde
D$ of $D$). Now  \[\frac{H^1(Y,{\bf C})}{F^1H^1(Y,{\bf
C})}\cong H^1(Y,{\cal O}_Y),\] and since $F^1$ yields an
exact functor on the category of mixed Hodge structures
(see \cite{D}), we get that \[\frac{H^1(Y,{\bf
C})}{F^1H^1(Y,{\bf C})}\cong\frac{H^1(U,{\bf
C})}{F^1H^1(U,{\bf C})}.\] Hence there is a natural map
\[{\rm Div}^0_S(Y)\to \frac{H^1(U,{\bf C})}{F^1H^1(U,{\bf
C})+H^1(Y,\Z (1))}\cong\frac{H^1(Y,{\bf C})}{F^1H^1(Y,{\bf
C})+H^1(Y,{\bf Z}(1))}\cong{\rm Pic}\,^0(Y).\] This map is
known to be just the natural map $D\mapsto{\cal O}_Y(D)$.
Hence \[\frac{H^1(U,{\bf C})}{F^1H^1(U,{\bf C})+H^1(U,{\bf
Z}(1))}\cong\frac{{\rm Pic}\,^0(Y)}{\mbox{subgroup
generated by ${\rm Div}^0_S(Y)$}}.\] Since $\limdir{S}{\rm
Div}\,^0_S(Y)$ is the group of all divisors on $Y$ which
are homologous to 0, which maps surjectively to ${\rm
Pic}\,^0(Y)$, the Lemma follows immediately.    \end{proof}
\B{rmk} Let $Y$ be a non-singular projective complex
variety. By the same argument in the proof of the Lemma,
one can see that \[\limdir{\mbox{$U\subset Y$}}
\frac{H^{2p-1}(U,{\bf C})}{F^pH^{2p-1}(U,{\bf C})+
H^{2p-1}(U,{\bf Z}(p))}\;=0\] (where now the limit is taken
over Zariski open subsets $U$ of $Y$ such that the
codimension of $Y-U$ is $\geq p$), if and only if
$J^p(Y)$ is algebraic \ie, the Abel-Jacobi map is
surjective.
\E{rmk}
\B{proof}$\!${\bf of the Theorem} Let $\cH^2(\Z(2))$
and $\cH^2/\cF^2$ be the Zariski sheaves on
$X$ associated to $U\mapsto  H^2(U,\Z(2))$ and
$U\mapsto  H^2(U,\C)/F^2$ respectively. Let
$\cH^{3}(\Z(2)_{\cD})$ denote the Deligne--Beilinson
cohomology Zariski sheaf on $X$ (\cf \cite{G} and
\cite{E}). We then have the following diagram of
sheaves
 $$\B{array}{ccccccccc}
0&\to &\cH^2(\Z(2))&\to&
H^2(\C (X),\Z(2)) &\to &\coprod_{x\in X^1}^{}
i_x H^1(\C (x),\Z(1))&\to&\cdots\\
&&\downarrow&&\downarrow&&\downarrow&&\\
0&\to &\cH^2/\cF^2&\to&
H^2(\C (X))/F^2 &\to &\coprod_{x\in X^1}^{}
i_x H^1(\C (x))/F^1&\to &0\\
&&\downarrow&&\downarrow&&\downarrow&&\\
0&\to &\cH^3(\Z(2)_{\cD})&\to&
H^3(\C (X),\Z(2)_{\cD})&\to &\coprod_{x\in X^1}^{}
i_x H^2(\C (x),\Z(1)_{\cD})&\to &0\\
&&\downarrow&&\downarrow&&\downarrow&&\\
&&0&&0&&0&&
\E{array}$$
obtained by comparing the arithmetic resolutions of the left hand column
(see \cite{BV} for that of $\cH^2/\cF^2$) via the standard long
exact sequences in Deligne-Beilinson cohomology (and taking into
account that singular cohomology of an affine vanishes in
degrees greater than its algebraic dimension). Now, because of the
Lemma, we have that  $$H^2(\C (x),\Z(1)_{\cD}) \df
\limdir{\mbox{$U\subset \overline{\{x\}}$}}
H^{2}(U,\Z(1)_{\cD})\cong  \limdir{\mbox{$U\subset
\overline{\{x\}}$}} \frac{H^1(U,{\bf C})}{F^1H^1(U,{\bf
C})+H^1(U,{\bf Z}(1))}\;=0.$$  Thus the right-bottom
corner in the diagram above vanishes and the sheaf
$\cH^{3}(\Z(2)_{\cD})$ is then identified by the constant
sheaf associated with the group $$H^{3}(\C
(X),\Z(2)_{\cD})\df\limdir{\mbox{$U\subset X$}}
H^{3}(U,\Z(2)_{\cD})\cong \limdir{\mbox{$U\subset X$}}
\frac{H^2(U,{\bf C})}{F^2H^2(U,{\bf C})+H^2(U,{\bf
Z}(2))}$$ where the limit is taken over all non-empty
Zariski open subsets of $X$.  Finally, if $p_g(X)=0$, it is
known (see \cite{BV} and \cite{E}) that
$$H^0(X,\cH^3(\Z(2)_{\cD}))\cong  \ker  (A_0(X)\by{\phi}
J^2(X)),$$ yielding the claimed identification.
 \E{proof}
\B{rmk} Note that for a surface $X$ with $p_g\neq 0$ we
have that $$\limdir{\mbox{$U\subset X$}}
\frac{H^2(U,{\bf C})}{F^2H^2(U,{\bf C})+H^2(U,{\bf
Z}(2))}\neq 0,$$ since this group maps
onto the ``Albanese kernel'', which is non-zero by
\cite{M}.

For a surface $X$ with $p_g=0$, the canonical
``residue map''  $$H^2(\C (X))/F^2\to \coprod_{x\in
X^1}^{}H^1(\C (x))/F^1$$ is injective, with quotient group
isomorphic to $H^3(X,\C)/F^2$ (see \cite{BV}); thus the
integral lifting of the ``residue'' of a given $\omega\in
H^2(\C (X))/F^2$ always exists, but, in order to lift
$\omega$ itself to a class in $H^2(\C(X),\Z(2))$, one
should know that the lifting of the ``residue'' of $\omega$
yields zero in the group of zero-cycles as well as in
$H^3(X,\Z(2))$.
 \E{rmk}

\vspace{1cm}

\B{flushright}
Dipartimento di Matematica\\
Universit\`a di Genova\\
Via L.B.Alberti, 4\\
16132 -- {\sc Genova}\\
{\it Italia}\\[8pt]

Tata Institute of Fundamental Research\\
School of Mathematics\\
Homi Bhabha Road\\
400 005 -- {\sc Bombay}\\
{\it India}
\E{flushright}

 \end{document}